\documentclass{article}
\usepackage{amsmath}
\usepackage{graphicx}
\def\bel{\begin{equation}\label}
\def\eeq{\end{equation}}
\newtheorem{Definition}{Definition}[section]
\newtheorem{Theorem}{Theorem}[section]
\newtheorem{Remark}{Remark}[section]
\newtheorem{Lemma}{Lemma}[section]
\newtheorem{Proposition}{Proposition}[section]
\newtheorem{Corollary}{Corollary}[section]
\newtheorem{Example}{Example}[section]

\title{The value function of an asymptotic exit-time optimal control problem}
\author{M.Motta and C.Sartori\\
Dipartimento di Matematica\\
Via Trieste, 63 - 35121 Padova, Italy\\
Telefax (39)(049) 8271428\\
e-mail: motta@math.unipd.it \\ sartori@math.unipd.it}
\def\fudge{\mathchoice{}{}{\mkern.5mu}{\mkern.8mu}}
\def\bbc#1#2{{\rm \mkern#2mu\vbar\mkern-#2mu#1}}
\def\bbb#1{{\rm I\mkern-3.5mu #1}}
\def\bba#1#2{{\rm #1\mkern-#2mu\fudge #1}}
\def\bb#1{{\count4=`#1 \advance\count4by-64 \ifcase\count4\or\bba A{11.5}\or
\bbb B\or\bbc C{5}\or\bbb D\or\bbb E\or\bbb F \or\bbc G{5}\or\bbb H\or
\bbb I\or\bbc J{3}\or\bbb K\or\bbb L \or\bbb M\or\bbb N\or\bbc O{5} \or
\bbb P\or\bbc Q{5}\or\bbb R\or\bbc S{4.2}\or\bba T{10.5}\or\bbc U{5}\or%
\bbb P\or\bbc Q{5}\or\bbb R\or\bba S{8}\or\bba T{10.5}\or\bbc U{5}\or
\bba V{12}\or\bba W{16.5}\or\bba X{11}\or\bba Y{11.7}\or\bba Z{7.5}\fi}}

\def \R {{\bb R}}
\def \C {{\mathcal  T}}
\def \N {{\bb N}}

\def \vv{\vskip 0.5 truecm}

\begin{large}
\begin{document}
\maketitle
\begin{abstract}{We consider a  class  of  exit--time   control problems for nonlinear systems with a nonnegative vanishing Lagrangian. In general,  the associated PDE may have multiple solutions, and  known regularity and  stability  properties  do not hold.   In this paper we  obtain  such properties and a uniqueness result  under some explicit sufficient conditions. We briefly investigate also the infinite horizon problem.}
\end{abstract}

\footnotetext  {$({\bf *})$ This research is partially supported by the Marie Curie ITN SADCO, FP7-PEOPLE-2010-ITN n. 264735-SADCO  and by the MIUR grant PRIN 2009 
"Metodi di viscosit\`a, geometrici e di controllo per modelli diffusivi nonlineari" (2009KNZ5FK)}

\section{Introduction}
Among the hypotheses under which the boundary value problem, (BVP),  
 \begin{equation}\label{Hin}
\left\{\begin{array}{c}{\mathcal H}(x,Du(x))= 0 \\
u=0\quad {\rm on}\,\,\partial \C,
\end{array}\right.
\end{equation}
with  \begin{equation}\label{IntH3}
 {\mathcal H}(x,u(x), D u(x))\doteq\sup_{a\in A}\{-\langle D u(x),f(x,a)\rangle-l(x,a) \}=0 
\end{equation}
and $\C$ a closed subset to $\R^n$ with compact boundary,  has a unique solution, the condition that $l\ge c_0>0$  together with small time local controllability, STLC,  around $\C$ plays a crucial role.
In this case,  under quite standard assumptions on the data,  the solution  to (\ref{Hin}) is represented as the value function  of an exit-time optimal control problem
with target  $\C\subset\R^n$,   trajectories governed by a   nonlinear control  system 
\begin{equation}\label{S}
\dot y(t)= f(y(t),\alpha(t)), \qquad y(0)=x  \quad(x\in\R^n)
\end{equation}
 and payoff given by
\begin{equation}\label{1.2}
{\mathcal J}(t,x,\alpha)=\int_0^tl(y(s),\alpha(s))\,ds,
\end{equation}
 where the control  $\alpha(t)$ belongs to the set $A\subset\R^m$  (assumed here to be compact).
More precisely,    for any $x\in\C^c\doteq\R^n\setminus{\mathcal T}$ the solution to (\ref{Hin}) is given by the  exit-time value function  
\begin{equation}\label{calvf}
  {\cal V}^f(x) \doteq\inf_{\{\alpha\in{\cal A}: \,\, t_x(\alpha)<+\infty\}}{\mathcal J}(t_{x}(\alpha), {x},\alpha) \quad(\le+\infty), 
 \end{equation}
where ${\cal A}$ is the set of measurable controls on $A$ and, for any $\alpha$, 
\begin{equation}\label{1.3}
t_x(\alpha) \doteq\inf\{t\ge0:      \ \  y_x(t,\alpha)\in{\mathcal  T}\} \quad(\le +\infty)  
\end{equation}
is the exit--time from $\C^c$. 

   In many interesting applications though, as for instance   the F\"uller or the shape from shading problems, we have $l\ge0$ and  the set
 \bel{Z}
{\mathcal{ Z}}\doteq\{x: \ l(x,a)=0 \ \text{for some } \,  a\in A\},
\eeq
 is non empty.  It is well known that in this case,  without additional hypotheses,  there is no hope to have a unique solution to (\ref{Hin}), even among the continuous, nonnegative functions. 
  
   In this paper we consider  an {\it exit--time} problem with $l$ nonnegative and  
   ${\mathcal T}\cap{\mathcal Z}\ne\emptyset$
and in Section \ref{infinite} we will see  that, in so doing, we can also cover some   {\it  undiscounted infinite horizon} problems like the  LQR  problem.
   Starting from the trivial remark that  we might have minimizing  trajectories approaching ${\mathcal{ T}}$ in infinite time,       
we introduce for any $x\in\C^c$  the following  {\it asymptotic} exit-time value function,
 \begin{equation}\label{calva}
 {\mathcal{ V}}({x})\doteq 
 \inf_{\alpha\in{\mathcal{ A}}({x})}{\mathcal J}(t_{x}(\alpha), {x},\alpha )  \quad(\le+\infty), 
\end{equation}
where 
$$
{\mathcal{ A}}({ x})\doteq\{\alpha\in{\mathcal A}:   \    \liminf_{t\to  t^-_{ x}(\alpha)} \text{dist}(y_{ x}(t,\alpha),\C)=0\}. 
$$
Then we  characterize ${\mathcal V}$ as the unique nonnegative  solution  to (\ref{Hin}), under  some global assumptions discussed below, and a  {\it special local asymptotic controllability} hypothesis  on the target,  also involving  the  Lagrangian   (see   \cite{MR} and (LACL) in Section \ref{comparison} below). Replacing the classical local small time controllability,  STLC,  on $\C$  by this  {\it weaker}  assumption  implies that,  while  ${\mathcal V}$ will be well defined in a neighborhood of the target,  in general
${\mathcal V}^f$ will not   be finite there.

\noindent Furthermore,  we perturb the (BVP) with more regular problems having a unique solution, show that their limit, say ${\mathcal{ U}}$,  not   coinciding in general with ${\mathcal{ V}}$,  can be represented as the value function of a constrained optimization problem.   Moreover, we give sufficient conditions for the equality ${\mathcal{ U}}\equiv{\mathcal{ V}}$, which essentially require the continuity of either ${\mathcal{ U}}$ or  ${\mathcal{ V}}$  on $\partial\C$.   
 
In more detail,  optimality principles obtained  in \cite{M} and  \cite{Sor2} tell us that ${\mathcal V}$ is the maximal solution to (\ref{Hin})
 when  it is continuous on $\partial\C$. The minimal  nonnegative  solution  to (\ref{Hin}) is represented by the value function:
\begin{equation}\label{calve}
 {\cal V}^m(x)\doteq \inf_{\alpha\in{\mathcal{ A}}}{\mathcal J}(t_{x}(\alpha), {x},\alpha ),
\end{equation}
where the minimization is done over trajectories not necessarily steering to the target.
Following this approach, we reduce the problem of uniqueness  essentially  to  the control theoretical questions of whether  ${\mathcal V}$ is continuous on $\partial\C$  and ${\mathcal V}\equiv {\mathcal V}^m$. 

Exploiting some recent results by Motta and Rampazzo \cite{MR}, we give suitable asymptotic controllability conditions  implying  the continuity  of   ${\mathcal U}$, introduced above,  and of ${\mathcal V}$,  on the target, but in general not in their whole  domains.  We also investigate  the global continuity  of ${\mathcal U}$ and ${\mathcal V}$  by introducing a kind of {\it turnpike condition},  that roughly states that trajectories not uniformly approaching the target, at least asymptotically,  are unaffordable (see e.g. \cite{Zas},  \cite{TZ}).  

 In order to have  ${\mathcal V}\equiv {\mathcal V}^m$, we introduce  some explicit sufficient  conditions   satisfied in many applications and generalizing several  previous hypotheses.   Let us remark that  many of the assumptions existing in the literature  ensuring uniqueness, {\it imply}   that   ${\mathcal V}\equiv {\mathcal V}^m$.  We refer to  Section 5 of \cite{M},  \cite{Ma} and  Remark \ref{Mal}  below  for an analysis  of some of them. 
 
Our uniqueness and stability results  improve  some previous  research  under several aspects.   Substituting  the STLC with the (LACL) and using  the  optimality principles instead of classical comparison theorems,  allows us to have the uniqueness of a continuous solution among all the nonnegative solutions to (BVP). Moreover,  owing to the  (LACL),  we  can also extend our results to  infinite horizon problems, where  admissible trajectories approach asymptotically  some set $\C\subset {\mathcal Z}$, noticeably the origin in the LQR problems.  When uniqueness fails, we give sufficient conditions in order to characterize ${\mathcal V}$ as  the limit of   (unique) nonnegative solutions of perturbed boundary value problems. We point out that  any stability result  has to be proved directly, as it cannot rely on the standard viscosity approach, based on the uniqueness of the solution to (\ref{Hin}). 

Without aiming to be exhaustive,  for the uniqueness issue and for an insight into   many  applications in which $l$ is not strictly positive, we refer  to   \cite{IR}, \cite {CSic},  \cite{Sor2},  \cite{M},  \cite{Ma}, \cite{CDP}, \cite{DL},  and \cite{G},  also concerning   unbounded controls and infinite horizon problems.    In particular, in  \cite{IR} and  \cite {CSic} special Hamiltonians are considered;  in \cite{Ma} just  the function ${\mathcal V}^f$ is characterized; in \cite{CDP} and  \cite{DL} uniqueness is obtained in the class of continuous functions  with bounded subdifferential, and  in \cite{G} among the convex functions. 

  \noindent   Finally,  let us mention that, in the  companion paper \cite{MS}, we extend the research begun here to the case of a non compact control set
and unbounded data.
 
The paper is organized as follows. In Section \ref{comparison} we characterize ${\mathcal V}$ as unique nonnegative solution of the (BVP).   Section \ref{App} is devoted to the approximation of ${\mathcal V}$. In Section \ref{Continuity results} we  give sufficient conditions  for the  continuity   either of  ${\mathcal V}$ in its domain,  or of the limit function  of the penalized problems, ${\mathcal{ U}}$, on $\partial\C$.  

 \noindent{\bf Notations.} Let $D\subset \R^N$ for some $N\in\N$.  $\forall r>0$  we 
 denote by $D_r$ the closed set $\overline{B(D,r)}$, while $D_r^c=\R^N\setminus D_r$.  $\overset{\circ} D$  is the interior of $D$. Moreover, $\chi_D$ denotes the characteristic function
of $D$, namely for any $x\in \R^N$ we set $\chi_D(x)=1$ if $x\in D$ and $\chi_D(x)=0$ if $x\notin D$.
For any function $u:\R^n\setminus\overset{\circ}{\C}\to\ R\cup\{+\infty\}$, we   denote the set $\{x\in  \R^n\setminus\overset{\circ}{\C}: \ u(x)<+\infty\}$ by $Dom(u)$. $[0,+\infty[\doteq\R_+$. A function $\omega:\R_+\times\R_+\to\R_+$
is called a {\it modulus} if: $\omega(\cdot,R)$ is increasing in a neighborhood of $0$, 
continuous at $0$, and $\omega(0, R)=0$ for every $R>0$; $\omega(r,\cdot)$ is increasing for every $r$.    Let $\Omega\supset\C$ be an open set and let  $U:\Omega\setminus\overset{\circ} \C \to\R_+$  be a locally Lipschitz  function. Then  $D^* U(x)\doteq  \{ p\in\R^n: \   p=\lim_{k}\nabla {U}(x_k), \   x_k\in diff(U)\setminus\{x\}, \ \lim_k x_k=x\}$ is  
the   {\it set of limiting gradients} of $U$ at $x$ (here $\nabla$ denotes the gradient operator and $diff(U)$ is the set of differentiability points of $U$).    For the notion of locally  semiconcave function and of viscosity solution  we refer e.g. to \cite{CS}, \cite{BCD}. ${\cal KL}$  denotes  the set of all continuous functions  $\beta:\R_+\times\R_+\to\R_+$ such that: (1)\, $\beta(0,t)=0$ and $\beta(\cdot,t)$ is strictly increasing  and unbounded for each $t\ge0$; (2)\, $\beta(r,\cdot)$ is decreasing  for each $r\ge0$; (3)\, $\beta(r,t)\to0$ as $t\to+\infty$ for each $r\ge0$.

\section{Uniqueness}  \label{comparison}
In this section we introduce  sufficient conditions,  under which we can   characterize ${\mathcal V}$ as unique nonnegative solution to (\ref{Hin}).  The present assumptions  generalize and in some sense unify several previous hypotheses, introduced either for exit-time or for undiscounted infinite horizon problems.  In particular, our uniqueness result does not require neither the STLC around the target nor ${\mathcal Z}\subset\C$  (see also Remark \ref{Mal} below).   We end the section with two simple, illustrative examples.

\vv
Let us begin by stating the hypotheses assumed throughout the whole  paper and the precise  definition of the boundary value problem. 

The control set $A\subset\R^m$ is compact and the target set ${\cal T}\subset\R^n$ is closed, with compact boundary. 
The function $l:\R^n\times A\to\R_+$ is continuous.  Moreover,   $f:\R^n\times A\to\R^n$  is continuous and there exist  $M>0$, and for any $R>0$,  there is some $L_R>0$   such that  
\begin{equation}\label{hp0}
\begin{array}{l}
|f(x_1,a)-f(x_2,a)|\le L_R\,|x_1-x_2|,\\
  |l(x_1,a)-l(x_2,a)|\le L_R\,|x_1-x_2| \quad\forall x_1, \ x_2\in\R^n, \ \forall a\in A, \\ 
|f(x,a)|\le M(1+|x|) \quad \forall x,  \ \forall a\in A.
\end{array}
\end{equation}
Hence for any $x\in\R^n$ and   for any measurable control $\alpha\in{\mathcal A}$,   (\ref{S}) admits just one solution,  defined on the whole interval $\R_+$. We use  
$y_{x}(\cdot,\alpha)$ (or, when no confusion may arise, $y_x(\cdot)$) to  denote such a  solution.

\begin{Definition}\label{DefBVP} {\bf (BVP)}  {\rm [M]} 
Any    function $u:\R^n\setminus \overset{\circ}{\mathcal{T}}\to\R\cup\{+\infty\}$ verifying   $u_*(x)\ge 0$ on $\partial\mathcal{T}$ and such that $u_*$ is a viscosity supersolution of 
\bel{Heq}
{\mathcal H}(x,Du(x))= 0 
\eeq
in $\R^n\setminus\C$,  is called a {\rm supersolution to} {\rm (BVP)}.  Any pair $(u,\Omega)$ where  $\Omega\supset{\mathcal T}$ is an open set  and  $u:\Omega\setminus \overset{\circ}{\mathcal{T}}\to\R$ is a locally bounded function verifying $u^*(x)\le 0$ on $\partial\mathcal{T}$ and such that $u^*$ is a viscosity subsolution of (\ref{Heq}) in 
$\Omega\setminus\mathcal{T}$, is called a {\rm subsolution to } {\rm (BVP)} $($in $\Omega)$. 

\noindent  Any pair $(u,\Omega)$, where $u:\R^n\setminus \overset{\circ}{\mathcal{T}}\to\R\cup\{+\infty\}$ and $\Omega$ is an open set,  $\Omega\supset{\mathcal T}$,   is called   a {\rm solution to } {\rm (BVP)} $($in $\Omega)$ if $u$ is a supersolution  and  $(u,\Omega)$ is a subsolution to  {\rm (BVP)}. 
 \end{Definition}

We recall the optimality principles and some related results obtained in \cite{M}  (see Thms. 2.1,  4.2, 4.3). We refer also to the works  \cite{Sor2}  and  \cite{Sor3},  where a strong formulation of the optimality principles has been first introduced and developed for this kind of problems.\footnote{We remark that in  \cite{M},    $l$ is also   satisfying
$ l(x,a)\le M(1+|x|)$ \,$\forall(x,a)\in\R^n\times A$,   for some $M>0$,  
but  this sublinear growth condition  can be removed,  as in  \cite{Sor3}.}  
   \begin{Proposition}  \label{EP0}   {\rm [M]}  Let $\mathcal{W} \in\{\mathcal{V},  \mathcal{V}^f,  \mathcal{V}^m\}$. 

(i) \, If  $\mathcal{W}$ is locally bounded in $Dom(\mathcal{W})$,  $Dom(\mathcal{W})$ is open and $\mathcal{W}^*\le0$ on $\partial{\mathcal T}$,  then $\mathcal{W}$ is a   subsolution to  {\rm  (BVP) } in $Dom(\mathcal{W})$.

(ii)  \, $\mathcal{W}$ is a  nonnegative   supersolution to  { \rm  (BVP)}. 
 \end{Proposition}

In the next statement we will use the following convexity hypothesis.
{\it \begin{itemize}
\item[{\bf (CV)}] For each  $x\in\R^n$,   the following set is convex:
\bel{convhp}
{\mathcal L}(x)\doteq \{ (\mu,\gamma)\in\R^{n+1}: \      \exists a\in A \ \text{s. t.}  \   \mu=f(x,a), \  l(x,a))\le\gamma\}. 
\eeq
\end{itemize}}
      
\begin{Proposition}  \label{EP3}     {\rm [M]} (i)  \, We have  ${\mathcal V}^m\le u$ for any  nonnegative and continuous   supersolution $u$  to { \rm  (BVP)}. If we assume {\rm (CV)}, then ${\mathcal V}^m$   is   l.s.c and it is  the minimal nonnegative   supersolution   to { \rm  (BVP)}.

(ii)  \,  If ${\mathcal V}$ is continuous on $\partial\C$, then  $({\mathcal V},\, Dom({\mathcal V}))$  is the maximal subsolution to {\rm (BVP)} among the pairs $(u,  Dom({\mathcal V}))$. \footnote{We recall that $Dom({\mathcal V})$ is an open set,  ${\mathcal V}$ is locally bounded and upper semicontinuous in view of Proposition \ref{UCS} below. }  
 \end{Proposition}
  \begin{Remark}\label{noCV}{\rm As usual,  if (CV) does not hold  the above result remains true if we replace ${\mathcal V}^m$  by the corresponding value function,  say ${\mathcal V}_r^m$,  obtained by taking the infimum over relaxed controls  (see \cite{M}). }
\end{Remark}

In \cite{Sor3}, there is a formally similar characterization of the maximal subsolution and minimal supersolution to (BVP) for a discontinuous Lagrangian. However,   in the  undiscounted  case considered here, those results are proved  when the  Lagrangian is  bounded below by a positive constant. 

 \vv
 We point out  that  (BVP) is a {\it free-boundary} value problem,  and that the exit-time value functions do not satisfy, in general,  the boundary condition
\bel{BCn}
\lim_{x\to\bar x}u(x)=+\infty \quad \forall \bar x\in\partial Dom(u).
\eeq
In the sequel, improving the results of \cite{M},   we  characterize  {\it the pair } $({\mathcal V}, Dom({\mathcal V}))$ as unique solution of (BVP), among the solutions verifying the boundary condition (\ref{BCn}). Disregarding such a restriction, we could still prove the uniqueness of ${\mathcal V}$    among the pairs $(u,\Omega)$ with $\Omega=Dom({\mathcal V})$.

\vv
Let us now state the following   continuity and uniqueness result,   whose proof follows from Theorem \ref{ET3} below. 
 This  general, but  quite theoretical   statement, is the   starting point for   handier results, given in the sequel.    
   
\begin{Theorem}\label{ET1}  Let ${\mathcal V}\equiv {\mathcal V}^m$.
  \begin{itemize}
\item[(i)] \, If  ${\mathcal V}$ is  continuous in $Dom({\mathcal V})$ and satisfies the boundary condition  {\rm (\ref{BCn})},  then  $\mathcal{V}$  is the unique nonnegative  viscosity solution to  {\rm({BVP})}  among the pairs $(u,\Omega)$, where $u$ is continuous in $\Omega$  and satisfies  {\rm (\ref{BCn})}.
\item[(ii)] \,  If {\rm (CV)} holds  and ${\mathcal V}$ is continuous on $\partial\C$, then  $(\mathcal{V}, Dom(\mathcal{V}))$  is  the unique nonnegative  viscosity solution to  {\rm({BVP})}  among the pairs $(u,\Omega)$, where $u$ satisfies {\rm (\ref{BCn})}. 
Moreover,  $\mathcal{V}$ is continuous.\footnote{Since ${\mathcal V}^{m}$ is lsc and ${\mathcal V}$ is usc, when $\mathcal{V}\equiv  {\mathcal V}^{m}$   (\ref{BCn}) is trivially satisfied.} 
\end{itemize}  
 \end{Theorem}
  One might wonder why the maximal subsolution of (BVP) is $\mathcal{V}$ and not $\mathcal{V}^f$, obviously larger.  This is not a contradiction, because,  if  $\mathcal{V}^f$  is  continuous on  $\partial\mathcal{T}$, so that it solves (BVP),  then $\mathcal{V}^f\equiv \mathcal{V}$ by Theorem \ref{P11} below. Let us stress however, that, using the above uniqueness result we can characterize the solution to (BVP) also in situations where   $\mathcal{V}<\mathcal{V}^f$, as shown by Examples \ref{Un1}, \ref{Un3} at the end of the section.

 \begin{Remark}\label{Mal}{\rm In the literature many uniqueness results for the solution of (BVP) concerning  the function  ${\mathcal V}^f$, are proved under hypotheses that imply  ${\mathcal V}\equiv{\mathcal V}^m$. For instance, this is easily seen in \cite{Ma1}, where, for all $x\in\C^c$ and $\alpha\in{\mathcal A}$, one supposes that
\vskip 0.2 truecm
$\int_0^{+\infty}l(y_x(t,\alpha),\alpha(t))\,dt<+\infty \quad  \Longrightarrow \quad \lim_{t\to+\infty}y_x(t,\alpha)\in\C.$
\vskip 0.2 truecm
\noindent 
Also the following hypotheses, used in \cite{Ma},
\vskip 0.2 truecm
 
a)  the trajectories that have a finite cost must stay in a bounded set,

b) $\forall t>0$ and $\forall \alpha\in{\mathcal A}$ one has $\int_0^tl(y_x(s,\alpha),\alpha(s))\,ds>0$,
\vskip 0.2 truecm
\noindent
 together with (CV), with a little bit of work can be shown to imply ${\mathcal V}\equiv{\mathcal V}^m$.  
 
 As discussed in Remark \ref{infH} in Section \ref{infinite}, in many undiscounted  infinite horizon problems, seen as asymptotic exit-time problems for  a suitable target (as the LQR problem),  condition ${\mathcal V}\equiv{\mathcal V}^m$ is naturally verified. 
 
 }\end{Remark}

Since   (BVP)  is a  free boundary   problem,  for any solution pair $(u, \Omega)$ we introduce  the Kruzkov transform   
$W(x)\doteq\Psi(u(x))\doteq1-{\rm e}^{-uv(x)}$, leading to another boundary value problem in $\R^n\setminus{\mathcal T}$,  whose solution, when unique,  simultaneously gives both  $u$  and $\Omega\doteq Dom(u)$.  
 
More precisely, the Hamiltonian associated to $W$ is
\begin{equation}
\label{HJB1}\begin{array}{l}
\mathcal{K}(x,u,p)
\doteq\sup_{a\in A}\{-\langle p,f(x,a)\rangle-l(x,a)+ l(x,a) u \} 
\end{array}
\end{equation}
and we consider the following  boundary value problem, in short, {\rm(BVP$\mathcal{K}$)},  
\begin{equation}\label{BVPK }
\left\{\begin{array}{l}
{\mathcal K}(x,W(x),DW(x))= 0\quad \qquad \text{in } \  \R^n\setminus\mathcal{T}\\
 W(x)= 0\qquad\qquad\quad\qquad \text{on } \ \partial\mathcal{T},
\end{array}\right.\end{equation}
where super- and subsolutions are defined analogously to Definition \ref{DefBVP}. 
 
\begin{Remark}\label{EPK} {\rm From \cite{M}, the statements of Propositions \ref{EP0}  and  \ref{EP3}  can be reformulated   in terms of the  Kruzkov transforms of the exit-time value functions, in the {\it whole} space  $\R^n\setminus{\mathcal T}$ (see  Cor. 3.1 and Thm. 4.3 in \cite{M}).  Incidentally,  these results hold without assuming the boundary condition (\ref{BCn}). Hence they  are not trivial, since,  given a subsolution $(u,\Omega)$ to (BVP),  $\Psi(u)$ is not  in general a subsolution to  {\rm(BVP$\mathcal{K}$)}  in $\R^n\setminus{\mathcal T}$ but just in $\Omega\setminus{\mathcal T}$. }
\end{Remark}  
   
We have the following uniqueness result in $\R^n$.    
\begin{Theorem}\label{ET3} Let ${\mathcal V}\equiv {\mathcal V}^m$. 
\begin{itemize}
\item[(i)] \, If  ${\mathcal V}$ is  continuous in $Dom({\mathcal V})$ and satisfies the boundary condition  {\rm (\ref{BCn})}, then there is a unique continuous, nonnegative viscosity solution $W$  to  {\rm(BVP$\mathcal{K}$)}. Moreover, 
${\mathcal V}\equiv\Psi^{-1}(W)=-\log(1-W)$ and $Dom(\mathcal{V})=\{x:\ \ W(x)<1\}$.
\item[(ii)] \,  If {\rm (CV)} holds  and ${\mathcal V}$ is continuous on $\partial\C$, then  there is a unique nonnegative viscosity solution $W$  to  {\rm(BVP$\mathcal{K}$)} which turns out to be  continuous. Moreover, 
${\mathcal V}\equiv\Psi^{-1}(W)=-\log(1-W)$ and $Dom(\mathcal{V})=\{x:\ \ W(x)<1\}$.
\end{itemize}  
\end{Theorem} 

\noindent  {\it Proof.}  We prove just (ii), the proof of  (i) being similar and actually simpler.  By Proposition  \ref {EP3}  and Remark  \ref{EPK}, for any solution $W$ to {\rm(BVP$\mathcal{K}$)} we get
 $$
 \Psi({\mathcal V}^m)(x)\le W_*(x)\le W(x)\le W^*(x)\le  \Psi({\mathcal V})(x) \quad \forall x\in\R^n\setminus\C.
 $$
Thesis (ii) follows now easily.   

 \vv
 Explicit sufficient conditions for the equality ${\mathcal V}={\mathcal V}^m$ are hypotheses (SC1), (SC2) below.

{\it \begin{itemize}
\item[{\bf (SC1)}]  There exists a Lyapunov function $U: \R^n\setminus\overset{\circ} \C  \to\R_+$,  $C^1$ in $ \R^n\setminus\overset{\circ}\C$, positive definite,  proper on $\C^c$ and  such that   $\forall x\in\C^c$,
\begin{equation}\label{C1s}
\sup_{ a\in A} \left\{   \langle  \nabla U(x), f(x,a)\rangle  \right\}\le -m({\bf d}(x))
\end{equation}
 for some continuous, increasing function  $m:]0,+\infty[\to]0,+\infty[$.
\end{itemize}}

{\it \begin{itemize}
\item[{\bf (SC2)}] there exists    a continuous, increasing function $c_2:]0,+\infty[\to]0,+\infty[$ such that 
\bel{H3'suf}
l(x,a)\ge  c_2({\bf d}(x))\   \qquad \forall (x,a)\in\C^c\times A.
\eeq 
\end{itemize}}

Condition (SC1) implies  that the control system (\ref{S}) is {\it  uniformly globally asymptotically stable, in short UGAS,  in $\C^c$,}  which roughly  means that {\it all}  trajectories with initial condition in a compact set,   approach $\C$ uniformly (see e.g.  \cite{BaRo} and the references therein).   We point out that (SC1) allows the Lagrangian to be zero outside the target. 

Hypothesis (SC2),  involving just the cost,   implies instead  that the set ${\cal Z}$ (see (\ref{Z}))  is a subset of $\C$.  It is  the simplest generalization of a strictly positive Lagrangian,  considered e.g. in \cite{Sor1},  in the framework of differential games.   Assumption (SC2) for $\C\equiv\{0\}$ is also  satisfied in LQR problems, where $l(x,a)=x^TQx+a^TRa$ and the matrices $Q$ and $R$ are  symmetric and positive definite. 

\begin{Proposition}\label{CorUn}   If   either   {\rm (SC1)} or  {\rm (SC2)} holds, then ${\mathcal V}={\mathcal V}^m$.
\end{Proposition}
 \noindent{\it Proof.}    It is immediate to see that condition  (SC2) yields the equality  ${\mathcal V}\equiv{\mathcal V}^m$, since ${\mathcal V}\equiv{\mathcal V}^m$ is equivalent to 
\bel{Cond=}
 \forall x\in\C^c: \quad \int_0^{+\infty}l(y_x(t),\alpha(t))\,dt \ge{\mathcal V}(x)  \quad\forall \alpha\in{\mathcal A}\setminus{\mathcal A}(x),
\eeq
and  (SC2) implies that
$$
 \int_0^{+\infty}l(y_x(t),\alpha(t))\,dt=+\infty  \quad\forall x\in\C^c, \ \forall \alpha\in{\mathcal A}\setminus{\mathcal A}(x).
 $$
Condition (SC1) instead, yields (\ref{Cond=}) because ${\mathcal A}\setminus{\mathcal A}(x)=\emptyset$ for all $x\in\C^c$.
 
\vv
 In order to state a special  local  asymptotic controllability condition,  in short (LACL),   introduced in \cite{MR} and sufficient to obtain the continuity of ${\mathcal V}$ on partial $\C$,  let us  recall the notion of (local) Minimum Restraint Function  from \cite{MR}.   For some terminology  borrowed from nonsmooth analysis we refer to the Notation.  Here   $h:\R^n\times A\to\R_+$ is an arbitrary   continuous Lagrangian.  
\begin{Definition}\label{CMRL}{\rm [MR]}   Given an open set $\Omega\subset\R^n$, $\Omega\supset\C$  we say that $U:\Omega\setminus\overset{\circ} \C \to\R_+$ is a {\rm local Minimum Restraint Function,}  in short, a local  MRF  for $h$,  if $U$  is  continuous on $\Omega\setminus\overset{\circ}\C$, locally  semiconcave, positive definite,  proper\,\footnote{$U$ is said {\it positive definite on $\Omega\setminus\C$} if  $U(x)>0$ \,$\forall x\in\Omega\setminus\C$ and $U(x)=0$ \,$\forall x\in\partial\C$.  $U$ is called {\it proper on $\Omega\setminus\C$} if $U^{-1}(K)$ is compact for every compact set  $K\subset\R_+$.}  on $\Omega\setminus\C$, $\exists U_0\in]0,+\infty]$ such that
 $$
\lim_{x\to x_0, \  x\in\Omega}U(x)=U_0 \ \  \forall x_0\in\partial\Omega;  \quad
U(x)<U_0 \quad\forall x\in\Omega\setminus\overset{\circ} \C,
$$
and, moreover,   $\exists{{k}}> 0$ such that, for every $x\in \Omega\setminus\C$,  
\begin{equation}\label{C1}
\min_{ a\in A} \left\{   \langle p, f(x,a)\rangle +k\,h(x,a) \right\}<0
 \qquad \forall p\in D^* U(x),
\end{equation}
   where  $D^* U(x)$ is   the    set of limiting gradients of $U$ at $x$. 
\end{Definition}
Let us observe that any MRF is a Control Lyapunov function for the system w.r.t. $\C$, which yields local asymptotic controllability to $\C$.
In the sequel, fixed a continuous  function $h\ge 0$,  we will often use the following hypothesis.
{\it \begin{itemize}
\item[{\bf (LACL)}]  There exists a local MRF  $U$ for $h$,   as introduced in Definition \ref{CMRL}.  
\end{itemize}}
By Theorem 1.1 in \cite{MR},  we have
\begin{Proposition}\label{CVT} If there exists a local   MRF $U$ for $l$,  then  
$
{\mathcal V}(x)\le U(x)/k
$
 in a neighborhood of the target. Hence ${\mathcal V}$ is   continuous  on $\partial\C$. 
\end{Proposition}
 Because of the degeneracy of $l$,  the continuity of  ${\mathcal V}$   on $\partial\C$ does not imply, in general,  the continuity in its   whole domain. Using a standard dynamic programming argument,  it is not difficult to prove  that in this case  ${\mathcal V} $ is upper semicontinuous.\begin{Proposition}\label{UCS}
Let  $\mathcal{V}$ be continuous on $\partial{\mathcal{ T}}$.  Then  $Dom(\mathcal{V} )$  is an open  set  and $\mathcal{V}$ is locally bounded and upper semicontinuous in it. 
\end{Proposition} 
   A general condition,   sufficient  for the propagation of the continuity of    ${\mathcal V}$,  will be given in Subsection \ref{RegRis}. In order to state the following explicit result,   let us  anticipate that both (SC1) and (SC2) imply such a propagation and the boundary condition {\rm (\ref{BCn})} (see Proposition \ref{ContVCS}).
 \begin{Corollary}\label{Unexpl} Assume  {\rm (CV)},   {\rm (LACL)}  for $l$ and either {\rm (SC1)} or {\rm (SC2)}. Then 
 \begin{itemize}
 \item[(i)]    there is a unique nonnegative viscosity solution $W$  to  {\rm(BVP$\mathcal{K}$)},   which turns out to be continuous.  Moreover, 
${\mathcal V}\equiv\Psi^{-1}(W)=-\log(1-W)$ and $Dom(\mathcal{V})=\{x:\ \ W(x)<1\}$;
 \item[(ii)] $(\mathcal{V}, Dom(\mathcal{V}))$  is  the unique nonnegative  viscosity solution to  {\rm({BVP})}  among the pairs $(u,\Omega)$, where $u$ satisfies {\rm (\ref{BCn})}.  Moreover,  $\mathcal{V}$ is continuous.
 \end{itemize} 
 If  {\rm (CV)} is not satisfied,  the above uniqueness results hold just among the continuous functions. 
 \end{Corollary} 
 
In the following  examples satisfying either (SC2) or (SC1), respectively,   we have uniqueness without the STLC on the target.
\begin{Example}\label{Un1}{\rm Let ${\mathcal T}\doteq\{0\}$,  consider  the scalar control system
$$ y'=-y\,a \quad\forall t>0, \qquad y(0)=x, \qquad a\in[0,1],
$$
 and define
$$
{\mathcal J}(t,x,\alpha)\doteq \int_0^{t}|y_x(t,\alpha)|\,dt.
$$
Since $y(t)=x\,e^{-\int_0^t|\alpha(s)|\,ds},$  for  every $x\ne0$ a  control $\alpha$  belongs to ${\mathcal A}(x)$ if and only if it verifies $\int_0^{+\infty}|\alpha(s)|\,ds=+\infty$. Hence ${\mathcal V}^f(x)=+\infty$,  while setting $\alpha\equiv1$ we can easily deduce that   
${\mathcal V}(x)=|x|$.

\noindent Notice that the system is asymptotically controllable to $\{0\}$,  $\{0\}$ is the unique zero of the Lagrangian and actually $l(x,a)=|x|$ satisfies (SC2). Moreover,  (CV) holds and  ${\mathcal V}$ is continuous everywhere. Incidentally, for any $k\in]0,1[$ the value function itself is a MRF for $l$. Therefore,  in view of Corollary \ref{Unexpl},    ${\mathcal V}(x)=|x|$ is the unique nonnegative solution of the boundary value problem 
$$
\max_{a\in[0,1]}\{ \langle Du, \, x\,a\rangle- |x|\}=0, \quad u(0)=0.
$$}
\end{Example}

 \begin{Example}\label{Un3}{\rm For any $x\in\R^2$ and any measurable function $\alpha:\R_+\to [0,1]$,  consider  the  control system 
$$ \left\{\begin{array}{l}y'(t)=-y(t)-y(t)\,\alpha(t) \quad\forall t>0, \\
y(0)=(y_1,y_2)(0)=( x_1, x_2)=x,
\end{array}\right.$$
and  the payoff
$$
 {\mathcal J}(t,x,\alpha)\doteq\int_0^{t}y_1^2(s,\alpha)\,ds \qquad\forall t>0,
$$
with the target $\C\doteq \{(0,0)\}$.
Since  $y_x(t,\alpha)=x\,e^{\int_0^t(-\alpha(s)-1)\,ds}$ in correspondence to any control $\alpha$, it is not difficult to prove that,   for every $x\ne0$,  ${\mathcal V}^f(x)=+\infty$ and 
$$
{\mathcal V}(x)\doteq \inf_{\alpha\in{\mathcal A(x)}}{\mathcal J}(t_x(\alpha),x,\alpha)=\frac{x_1^2}{4}. 
$$
Notice that the zero level set of the Lagrangian $l(x,a)=x_1^2$ is given by the unbounded set ${\mathcal Z}=\{(0,x_2):\  x_2\in\R\}$, but the control system is UGAS in $\R^2\setminus \{(0,0)\}$, so that condition (SC1) is verified. Moreover, (CV) holds  and it is easy to see that  the following function 
$$U(x)\doteq
\frac{x^2_1+x_2^2}{4}$$
is a MRF function for $l(x,a)=x_1^2$ that  verifies  (\ref{C1})
 for any $k\in]0,1[$.  

\noindent Hence in view of Corollary \ref{Unexpl},  ${\mathcal V}(x)=x_1^2/4$    is the unique nonnegative viscosity solution of the following (BVP): 
 $$\max_{a\in[0,1]}\{-\langle D u(x),-x-xa\rangle-x_1^2 \}=0 \quad\forall  x\in\R^2\setminus\{(0,0)\},\qquad u(0,0)=0.
 $$ 
}
\end{Example}

When  ${\mathcal V}\ne {\mathcal V}^m$,  owing to Proposition \ref{EP3}   we can still characterize ${\mathcal V}$  as maximal subsolution of (BVP).  Following an alternative approach,  we can obtain  ${\mathcal V}$  as limit of perturbed  boundary value problems (having a unique  solution). This is the topic of the next section. 
 
  \section{Approximation results} \label{App}
This section is mainly devoted to  discuss  whether and when the  exit-time problem (\ref{calva}) can be   approximated by  some perturbed   problems, giving also  {\it uniform} convergence conditions.  These results can be seen as  stability properties  for the (BVP), in the sense that   ${\mathcal V}$ can be selected as the solution which is the  limit of  value functions,  themselves characterized as unique solution of  suitable boundary value problems. 
 
\subsection{Penalized problems}
 Let  $\rho: \R^n\times A \to \R_+$ be a continuous function. Fix $\varepsilon>0$.  For every $x\in{\mathcal T}^c$, $\alpha\in{\mathcal A}(x)$,  let us define the  {\it $\varepsilon$-penalized  payoff,} 
 \begin{equation}\label{Je}
 {\mathcal J}^\rho_\varepsilon(t,x,\alpha)\doteq\int_0^{t}\left[l( y_{ x}(\tau,\alpha),\alpha(\tau))+\varepsilon\rho( y_{ x}(\tau,\alpha),\alpha(\tau))\right]\,d\tau
 \end{equation}
 and the corresponding  {\it $\varepsilon$-penalized value function,} 
 \begin{equation}\label{Vpen}
 {\mathcal V}^\rho_\varepsilon (x)\doteq  \inf_{\alpha\in{\mathcal{ A}}({x})} {\mathcal J}_\varepsilon(t_x(\alpha), x,\alpha). 
 \end{equation}
 We introduce also  the  {\it penalized value function}
\bel{Vrho}
{\mathcal U}^\rho(x)\doteq \inf_{\alpha\in{\mathcal{ A}_\rho}({x})} {\mathcal J}(t_x(\alpha), x,\alpha)\quad(\le+\infty),
\eeq
where 
 \bel{ro}
{\mathcal A}_\rho(x)\doteq\left\{\alpha\in{\mathcal A}(x): \ \ \int_0^{t_x(\alpha)} \rho(y_x(t,\alpha),\alpha(t))\, dt <+\infty \right\},
\eeq
which will play a crucial role in the sequel.

\noindent  Clearly, ${\mathcal V}(x)\le {\mathcal U}^\rho(x)\le {\mathcal V}^\rho_\varepsilon(x)$ and the inequalities may be strict,  as one can easily see in  Example \ref{Un1} choosing $\rho(x,a)\doteq |a|$ for every $(x,a)$.

For any $\varepsilon>0$, let ${\mathcal K}_\varepsilon(x,u,p)$ denote the Hamiltonian defined as ${\mathcal K}$ in (\ref{HJB1}) with $l$ replaced by $l+\varepsilon\rho$. 
 In view of Theorem \ref{ET3} and  Corollary \ref{Unexpl},   Theorem \ref{UnConv}  below implies the following stability result. 
\begin{Theorem}\label{Stab}   Let   $\rho: \R^n\times A \to \R_+$ be  a  continuous  function. Assume  {\rm (CV)}, {\rm (LACL)} for $l+\rho$\,\footnote{The last hypothesis can be replaced by  the weaker explicit sufficient conditions in  Theorem \ref{FTFF} below,  implying the  continuity of the limit function ${\mathcal U}^\rho$ on $\partial\C$. 
} and either  {\rm (SC1)} or {\rm (SC2)} for $l$ replaced by $l+\rho$. Then  for any $\varepsilon\in]0,1]$ there exists a unique nonnegative solution $W_\varepsilon$ to
\bel{BVPKe}
\left\{\begin{array}{l}
{\mathcal K}_\varepsilon(x,W(x),DW(x)) = 0\quad \qquad \text{in } \  \R^n\setminus\mathcal{T}\\
 W(x)= 0\qquad\qquad\quad\qquad \text{on } \ \partial\mathcal{T}.
\end{array}\right.
\eeq
Moreover,  as $\varepsilon\to0^+$ the $W_\varepsilon$ converge  to a function $W$ such that   ${\mathcal V}\equiv -\log(1-W)$ and $Dom(\mathcal{V})=\{x:\ \ W(x)<1\}$. If $W$ is continuous in $\R^n\setminus\overset{\circ}{\mathcal{T}}$, then the convergence is locally uniform. 

\noindent If  {\rm (CV)} is not satisfied,  for any $\varepsilon\in]0,1]$,   the function $W_\varepsilon$   is the unique solution to {\rm (\ref{BVPKe})}  just among the continuous functions.  
\end{Theorem} 

Notice that, if (SC1) holds,  owing to Proposition \ref{ContVCS},  ${\mathcal V}$ is continuous in its domain.  When instead  {\rm (SC2)} for   $l+\rho$ (not implying, in general, (SC2) for $l$) is assumed,   the global continuity of ${\mathcal V}$  is not guaranteed a priori.

\vv
 Choosing, e.g.,   $\rho(x,a)\doteq {\bf d}^r(x)$ for some integer $r\ge 1$, then for any $l\ge0$, $l(x,a)+\rho(x,a)$ satisfies (SC2) and, as  in  the following example, we can have uniqueness of the solution for the perturbed problems.  Let us remark that   in the next example the   trivial choice $\rho\equiv 1$ does not give an approximation of ${\mathcal V}$ (see also Proposition \ref{FFFT}). 
 
 \begin{Example}\label{UneV}{\rm Let ${\mathcal T}\doteq\{0\}$,  consider  the scalar control system of Example \ref{Un1},
$$ y'(t)=-y(t)\,\alpha(t) \quad\forall t>0, \qquad y(0)=x, \qquad \alpha(t)\in[0,1],
$$
 and define
$$
{\mathcal J}(t,x,\alpha)\doteq \int_0^{t}|y^2_x(s,\alpha)-y_x(s,\alpha)|\,ds \qquad \forall t>0.
$$
Notice that ${\mathcal V}^m\ne {\mathcal V}$. Indeed, implementing the control $\alpha\equiv 0$, one gets e.g. ${\mathcal V}^m(1)=0,$   while 
${\mathcal V}(1)=1/2$, being 
\bel{Ves}
{\mathcal V}(x)=\left|x-\frac{x^2}{2}\right|\qquad \forall x\in\R,
\eeq
 as we show below.  Let us  introduce for every $\varepsilon>0$   the $\varepsilon$-penalized value function
$$
{\mathcal V}^\rho_\varepsilon(x)\doteq\inf_{\alpha\in{\mathcal A}(x)}\int_0^{t_x(\alpha)}\left(|y^2_x(t,\alpha)-y_x(t,\alpha)|+\varepsilon|y_x(t,\alpha)|\right)\,dt.
$$
Setting $f(x,a)\doteq -xa$, $l(x,a)\doteq|x^2-x|$ and $\rho(x,a)\doteq |x|$,  in view of Theorem \ref{Stab},  for any $\varepsilon>0$  the Kruzkov transform  $W_\varepsilon \doteq 1-e^{-{\mathcal V}^\rho_\varepsilon }$ is the unique nonnegative solution to (\ref{BVPKe}) among the continuous functions, since (CV)is not verified. Moreover,  as $\varepsilon\to0^+$ the $W_\varepsilon$ converge   to  $W(x)=1-e^{\left|x-\frac{x^2}{2}\right|}$ for all $x\in\R$. Notice that the convergence is locally uniform in view of the continuity of $W$. Finally, ${\mathcal V}\equiv -\log(1-W)$. Indeed,  straightforward calculations show that  the  function
$$
U(x)\doteq\left|x-\frac{x^2}{2}\right|+|x|\qquad \forall x\in\R
$$
is a MRF function for $l +\rho$,  that  verifies (\ref{C1}) for any $0<k<1$
and  the lagrangian $|x^2-x|+|x|$ obviously verifies (SC2).  At this point it is easy to
show that
$$
{\mathcal V}^\rho_\varepsilon(x)=\left|x-\frac{x^2}{2}\right|+\varepsilon|x|\qquad \forall x\in\R.
$$
Finally,  ${\mathcal V}$ is given by (\ref{Ves}) and can be characterized as the  locally uniform limit of ${\mathcal V}^\rho_\varepsilon$.
}\end{Example}
  
 \vv
In Theorem \ref{UnConv} below we shall  give a representation formula for  the limit   of the penalized problems in terms of the value function ${\mathcal U}^\rho$ defined in (\ref{Vrho}).  Since  the zero level set of $l$ is arbitrary  such a limit  does not coincide  in general with     ${\mathcal V}$. In this case, problem (\ref{calva}) is sometimes said to exhibit the {\it Lavrentiev phenomenon} or, more precisely,  to have  a ${\mathcal A}(x)-{\mathcal A}_{\rho}(x)$ {\it Lavrentiev gap.}  As a first result, the next theorem shows that, on the one hand, such a phenomenon cannot occur   whenever  ${\mathcal V}\equiv{\mathcal V}^f$. On the other hand, it  provides a sufficient condition to avoid the  Lavrentiev gap in case  ${\mathcal V}\ne{\mathcal V}^f$.  Theorem \ref{P11} includes a result in this sense due to  Guerra and Sarychev,   \cite{GS},  concerning  the special case of  an affine system and $l(x,a)\equiv x'Px$, $P$ a symmetric and positive matrix,  under a  local stabilizability assumption, involving bounded controls. 
 
\begin{Theorem} \label{P11}
 \begin{itemize}
\item[(i)]  Let $x\in \C^c$. If  ${\mathcal V}(x)={\mathcal V}^f(x)$, then ${\mathcal U}^\rho(x)={\mathcal V}(x)$ for any continuous function $\rho$. 
\item[(ii)] If,  for some $\rho$, ${\mathcal U}^\rho$  is  continuous on $\partial{\mathcal T}$  then ${\mathcal U}^\rho\equiv {\mathcal V}$. Therefore, in particular, if  ${\mathcal V}^f$ is continuous  on $\partial\C$, then  
 ${\mathcal V}^f\equiv{\mathcal U}^\rho\equiv{\mathcal V}$ for any continuous, nonnegative function $\rho$.
\end{itemize}
\end{Theorem}

\noindent{\it Proof.}   The inequality ${\mathcal V}(x)\le {\mathcal U}^\rho(x)$ is obvious and,  if ${\mathcal V}(x)=+\infty$,   it implies immediately   ${\mathcal V}(x)= {\mathcal U}^\rho(x)$.  Let  $x\in{\mathcal T}^c$ with ${\mathcal V}(x)<+\infty$ and  let  $\eta>0$. 
 By the definition of ${\mathcal V}$, there exists a control $\tilde\alpha\in{\mathcal A}(x)$ such that
\begin{equation}\label{c1n}
\int_0^{t_x(\tilde\alpha)}l(y_x(t,\tilde\alpha),\tilde\alpha(t))\,dt< \mathcal{V} (x)+\eta.
\end{equation}

{\it Case 1\,}: ${\mathcal V}(x)={\mathcal V}^f(x)$, so that  for any $\eta>0$ we can assume $t_x(\tilde\alpha)<+\infty$.    By standard estimates there exists some $R>0$ such that $|y_x(t,\tilde \alpha)|\le R$ for all $t\in[0,t_x(\tilde\alpha)]$. 
Hence,    by continuity there is some $\bar M_R>0$  such that $\sup_{t\in[0, t_x(\tilde\alpha)]}\rho(y_x(t,\tilde \alpha), \tilde\alpha(t))\le \bar M_R$    and
$$
\int_0^{t_x(\tilde\alpha)}\rho(y_x(t,\tilde\alpha),\tilde\alpha(t))\,dt\le \bar M_R\, t_x(\tilde\alpha)<+\infty,
$$
Therefore  $\tilde\alpha\in\mathcal{A}_\rho(x)$ and  by (\ref{c1n}) we get that  ${\mathcal U}^\rho(x)\le  \mathcal{V} (x)+\eta$.    By the arbitrariness of $\eta>0$, this concludes the proof of  statement (i). 
 
 {\it Case  2\,}: ${\mathcal V}(x)<{\mathcal V}^f(x)$, so that $t_x(\tilde\alpha)=+\infty$ if $\eta<{\mathcal V}^f(x)-{\mathcal V}(x)$.
By the continuity of ${\mathcal{U}}^\rho$ on the compact set $\partial\mathcal{T}$,   there is  some  $\delta>0$ such that 
\bel{pA02}
{\mathcal{U}}^\rho(\bar x)<\eta \qquad \forall \bar x\in{\mathcal T}^c \quad\text{with}\quad {\bf d}(\bar x)<\delta.
\eeq
  Moreover for $\tilde\alpha\in{\mathcal A}(x)$  satisfying (\ref{c1n}), there is some $\bar  t< t_x(\tilde\alpha)$ such that
\begin{equation}\label{c11}
{\bf d}( y_x(\bar t,\tilde\alpha))< \delta/2.
\end{equation}
At this point, we get that
$\int_0^{\bar t}\rho(y_x(t,\tilde\alpha),\tilde\alpha(t))\,dt< +\infty$,  arguing as in  Case 1.
  Let $\tilde x\doteq y_x(\bar t,\tilde\alpha)$ and let  $\bar\alpha\in\mathcal{A}_\rho(\tilde x) $ be a control  such that 
$$
\int_0^{t_{ \tilde x}(\bar\alpha)}l(y_{ \tilde x}(t,\bar\alpha),\bar\alpha(t))\,dt<\eta,
$$
 which exists in view of  (\ref{pA02}). Then,  the control $\alpha(t)=\tilde\alpha(t)\chi_{[0,\bar t]}(t)+\bar\alpha(t-\bar t)\chi_{]\bar t,+\infty]}(t)$ 
 belongs to $\mathcal{A}_\rho(x)$  and  
\begin{equation}\label{c15}
\begin{array}{l}
\int_0^{t_x(\alpha)}l(y_x(t,\alpha),\alpha(t))\,dt  < \mathcal{V}(x)+3\eta,
\end{array}
\end{equation}
so that  ${\mathcal U}^\rho(x)<  {\mathcal{V}}( x)+3\eta$.  Statement (ii) is thus proved,   by the arbitrariness of $\eta>0$.

%
%
%
%
\vv
 For the proof of the uniform convergence result in  Theorem \ref{UnConv} below,  we need the following proposition, whose quite easy proof is omitted for the sake of brevity.
\begin{Proposition}\label{ConT}  Let   $\rho: \R^n\times A \to \R_+$ be  a  continuous  function  and let ${\mathcal V}^\rho_{\bar\varepsilon}$ be continuous on $\partial\C$ for some $\bar\varepsilon>0$. Then ${\mathcal U}^\rho$ is continuous on $\partial\C$;  for any $\varepsilon>0$,  ${\mathcal V}^\rho_{\varepsilon}$ is continuous on $\partial\C$;  $Dom(\mathcal{V}^\rho_{\varepsilon})=Dom({\mathcal U}^\rho)$ and it   is an open  set;    $\mathcal{V}^\rho_{\varepsilon}$ is locally bounded and upper semicontinuous in $Dom({\mathcal U}^\rho)$.  
\end{Proposition}

\begin{Theorem} \label{UnConv} Let   $\rho: \R^n\times A \to \R_+$ be  a  continuous  function.  
\begin{itemize}
\item[(i)] \,
One has
\bel{A01}
 \lim_{\varepsilon\to 0^+}{\mathcal V}^\rho_\varepsilon (x)= {\mathcal U}^\rho(x) \qquad \forall x\in{\mathcal T}^c
 \eeq
and, by Theorem \ref{P11},  ${\mathcal U}^\rho\equiv {\mathcal V}$ 
if either ${\mathcal V}^f= {\mathcal V}$ or ${\mathcal U}^\rho$ is continuous on $\partial\C$.
\item[(ii)] \,  If  there is some $\bar\varepsilon>0$ such that  ${\mathcal V}^\rho_{\bar\varepsilon}$ is continuous on $\partial\C$ and   ${\mathcal V}$ is continuous on its whole domain, then ${\mathcal U}^\rho\equiv {\mathcal V}$  and the  limit {\rm (\ref{A01})} is uniform on any compact set $Q\subset Dom ({\mathcal V})$. 
\end{itemize}
\end{Theorem}
 
\noindent{\it Proof.}    Let $\varepsilon>0$ and fix  $x\in{\mathcal T}^c$.  Since  for any $\alpha\in{\mathcal A}(x)$,    ${\mathcal J}_\varepsilon(t_x(\alpha), x,\alpha)<+\infty$ implies $\alpha\in \mathcal{ A}_\rho({x})$, then ${\mathcal U}^\rho(x)\le\mathcal{V}^\rho_\varepsilon(x)$ and in order to prove (i) it remains only to show that the strict inequality ${\mathcal U}^\rho(x)<\inf_{\varepsilon>0}\mathcal{V}^\rho_\varepsilon(x)$ cannot hold.  If ${\mathcal U}^\rho(x)=+\infty$,  the equality is obvious. If instead  ${\mathcal U}^\rho(x)<+\infty$, let us assume by contradiction that, for some $\eta>0$,  
$$
{\mathcal U}^\rho(x)< {\mathcal V}^\rho_\varepsilon (x)-3\eta
$$
for any   $\varepsilon>0$.  By the definition of  ${\mathcal U}^\rho$, there is some $\alpha\in{\mathcal{A}}_{\rho}(x)$     such that  
$$\int_0^{t_x(\alpha)}l(y_x(t,\alpha),\alpha(t))\,dt< {\mathcal U}^\rho(x) +\eta$$
and  $\bar\rho\doteq\int_0^{t_x(\alpha)}\rho(y_x(t,\alpha),\alpha(t))\,dt<+\infty$. Then one has 
\bel{veps}
  {\mathcal V}^\rho_\varepsilon (x)< \int_0^{t_x(\alpha)} l(y_x(t,\alpha),\alpha(t))\,dt+\eta
\eeq
  for any $\varepsilon\le \eta/\bar\rho$. The proof of (\ref{A01}) is thus concluded, since we immediately obtain the contradiction 
  $ {\mathcal U}^\rho(x)<   {\mathcal U}^\rho(x)-\eta.$  
  
  Let us now prove (ii).  In view of Proposition \ref{ConT}, for any $\varepsilon>0$, ${\mathcal V}^\rho_{\varepsilon}$  is  locally bounded and upper semicontinuous in its domain, which is an open set and coincides with $Dom({\mathcal U}^\rho)$. Moreover,  ${\mathcal U}^\rho\equiv {\mathcal V}$ in view of Theorem \ref{P11}.
Therefore,   $\{{\mathcal V}^\rho_{\varepsilon}-{\mathcal V}\}_{\varepsilon>0}$ is   a decreasing sequence of nonnegative and  upper semicontinuous functions defined in the open set $Dom({\mathcal V})$ and  converging to the null  function as $\varepsilon\to0^+$.
Now, the uniform convergence of the 
${\mathcal V}^\rho_{\varepsilon}$ to ${\mathcal V}$ on any compact subset $Q\subset Dom({\mathcal V})$ follows from an easy adaptation of   Dini's Theorem for continuous functions to  the upper semicontinuous case.
 

\vv
We conclude this subsection pointing  out that  in the special case $\rho\equiv 1$ the limit function  ${\mathcal U}^\rho$ of  $\varepsilon$-penalized problems  coincides  with the limit of the  so-called {\it $T$-finite time value functions,} defined  for any $T>0$, as follows:
 $$
{\mathcal V}_T(x)=\inf_{\{\alpha\in{\mathcal A(x)}: \ t_x(\alpha)\le  T\}}{\mathcal J}(t_x(\alpha),x,\alpha) \qquad\forall x\in\C^c.
$$
\begin{Proposition}\label{FFFT}  For every $x\in\C^c$,  
$$
{\mathcal V}^f(x)=\lim_{T\to+\infty}{\mathcal V}_T(x).
$$
 If ${\mathcal V}^f$  is continuous on $\partial\C$, then the above limit coincides with ${\mathcal V}$.  When in addition  ${\mathcal V}_T$ for some $T>0$ is continuous on $\partial\C$ and ${\mathcal V}$ is continuous in its domain, then the   above convergence is locally uniform.   
\end{Proposition}

\noindent{\it Proof.}   Fix  $x\in{\mathcal T}^c$.  Clearly,   ${\mathcal V}^f(x)\le\mathcal{V}_T(x)$ for any $T>0$ and we have just to prove  that  ${\mathcal V}^f(x)<\inf_{T>0}\mathcal{V}_T(x)$ cannot hold.  If ${\mathcal V}^f(x)=+\infty$,  the equality is trivial. If instead  ${\mathcal V}^f(x)<+\infty$, for any $\eta>0$ there is some  $\alpha\in{\mathcal{A}}(x)$     such that $T_\eta\doteq t_x(\alpha)<+\infty$ and  
$$
{\mathcal V}_{T_\eta}(x)\le\int_0^{t_x(\alpha)}l(y_x(t,\alpha),\alpha(t))\,dt\le {\mathcal V}^f(x) +\eta.
$$
This shows that the  limit holds.   The remaining result is a consequence of   Theorem \ref{UnConv}.

\subsection{Target approximations}
One could also consider an approximation {\it from below} of   ${\mathcal V}$, by fattening the target. For every $\delta>0$ and $x\in{\mathcal T}_\delta^c$, let us define   
$$t^\delta_x(\alpha)\doteq\inf\{t\ge0:\quad y_x(t,\alpha)\in{\mathcal T}_\delta\} \quad(\le+\infty)
$$
 and the ${\mathcal T}_\delta$-problem
$$
{\mathcal V}_{{\mathcal T}_\delta}(x)\doteq\inf_{\{\alpha\in{\mathcal A}:\quad t^\delta_x(\alpha)<+\infty\}}{\mathcal J}(t^\delta_x(\alpha),x,\alpha)\qquad(\le+\infty),
$$
where only trajectories reaching the   target  ${\mathcal T}_\delta$ {\it in finite time} are allowed. 

In the next proposition we show that ${\mathcal V}$ is  the natural limit of the ${\mathcal V}_{{\mathcal T}_\delta}$ as $\delta\to0^+$, as soon as 
it is continuous on $\partial\C$.  By Proposition \ref{CVT}, a sufficient condition for such a continuity is the (LACL) for $l$, which easily yields the  STLC for ${\mathcal T}_\delta$ $\forall\delta>0$   small enough,  but {\it not}  in general for ${\mathcal T}$.

\begin{Proposition}\label{targetA}  For every $x\in\C^c$,  assuming ${\mathcal V}$ continuous on $\partial\C$, we have
$$
{\mathcal V}(x)=\lim_{\delta\to 0}{\mathcal V}_{{\mathcal T}_\delta}(x)
$$
and this convergence is  uniform in $Dom({\mathcal V})$.   
\end{Proposition}

\noindent {\it Proof. }    Fix $\eta>0$. The continuity of $ \mathcal{V}$ on the compact set $\partial\mathcal{T}$ implies that there is some $\bar\delta>0$ such that $\forall \bar x\in {\mathcal T}_{\bar\delta}\setminus{\mathcal T}$ one has   $\int_0^{t_{\bar x}(\bar\alpha)}l(y_{\bar x}(t,\bar\alpha),\bar\alpha(t))\,dt\le\eta/2$ for some $\bar\alpha\in {\mathcal A}(\bar x)$. Let $S(x)\doteq\sup_{\delta>0}{\mathcal V}_{{\mathcal T}_\delta}(x).$ Obviously $S(x)\le\mathcal{V}(x),$ for all $x\in \mathcal{T}^c,$ and
if $S(x)=+\infty$ then $S(x)=\mathcal{V}(x).$ Therefore let $S(x)<+\infty$ and for every
$\delta\in]0,\bar\delta[$,  choose $\alpha_{\delta}\in\mathcal{A}$ satisfying
\begin{equation}\label{c1}
\int_0^{ t^\delta_x(\alpha_{\delta})}l(y_x(t,\alpha_{\delta}),\alpha_{\delta}(t))\,dt\le \mathcal{V}_{{\mathcal T}_\delta}(x)+\eta/2.
\end{equation}
Set $\bar t\doteq t^\delta_x(\alpha_{\delta})$,  $\bar x\doteq y_x(\bar t,\alpha_\delta)\in\mathcal{T}_\delta$ and  
 $\alpha(t)\doteq\alpha_\delta(t)\chi_{[0,\bar t]}+\bar \alpha(t-\bar t)\chi_{[\bar t,+\infty[}.$ Then $\alpha\in\mathcal {A}(x)$ and
 \begin{equation}\label{c2}\begin{array}{l}
\mathcal{V}(x)\le\int_0^{t_x(\alpha)}l(y_x(t,\alpha),\alpha(t))\,dt\le   \int_0^{  \bar t}l(y_x(t,\alpha_\delta),\alpha_\delta(t))\,dt+\eta/2 \\
\qquad\qquad\qquad\qquad\qquad\qquad\qquad\qquad\le \mathcal{V}_{{\mathcal T}_\delta}(x)+ \eta\le S(x)+  \eta\end{array}
\end{equation}
and, being $\eta>0 $ arbitrary,  the equality is proved. As we see from the proof, this convergence is uniform in $Dom ({\mathcal V})$ with no further assumptions.

\vv
It is easy to see that sufficient conditions in order to have uniqueness to
(BVPK) (see e.g. Corollary \ref{Unexpl}) yield the uniqueness to the boundary value problem naturally associated to the Kruzkov transform of ${\mathcal V}_{{\mathcal T}_\delta}$, for any $\delta>0$ sufficiently small.

\section{Continuity results}\label{Continuity results}

In the present section, we give sufficient condition for the continuity of ${\mathcal V}$ in its whole domain.  Furthermore, exploiting some results of \cite{MR}, we give some sufficient conditions in order to have the continuity of ${\mathcal U}^\rho$ on $\partial{\mathcal T}$ and we particularize  the results for   the value function ${\mathcal V}^f$.  

\subsection{Global continuity}\label{RegRis}
 The continuity of  the function   ${\mathcal V}$  on the target does not propagate in general to the whole domain,  but,  as stated in   Proposition \ref{UCS}, it   implies just  the upper semicontinuity of  ${\mathcal V}$.     In Theorem \ref{ContV}  we  prove that ${\mathcal V}$ is  globally continuous assuming its continuity  on the target and the {\it turnpike-type condition} (TPC)  below. Turnpike conditions have their roots in economic growth theory.   More recently their use has been extended to  a wide range of variational and optimal control problems. The novelty here is to introduce such a kind of notion  in order to  study the continuity of the value function ${\mathcal V}.$  We refer to the book \cite{Zas}  and to \cite{TZ} for  interesting surveys on the subject.
 
 \noindent 
 Loosely speaking  (TPC)  says that for any fixed $\delta$-neighborhood of the target,  for every  $x\in Dom({\mathcal V})$  one can select a nearly optimal asymptotic trajectory reaching  $\C_\delta$ in finite time $T$,    depending just on  ${\mathcal V}(x)$ and $\delta$ (that is, uniformly w.r.t. $x$ and   $\alpha$).  
 
{\it \begin{itemize}
\item[{\bf (TPC)}]     $\forall R$, $\eta$, $\delta>0$,   there exists some  increasing function $T(\cdot):\R_+\to\R_+$ such that for every  $x\in\C^c_\delta\cap  Dom({\mathcal V})$,  ${\bf d}(x)\le R$, 
there is a control  $\alpha\in{\mathcal A}(x)$ verifying
\bel{TK}
\begin{array}{l} 
\int_0^{t_x(\alpha)}l(y_x(t,\alpha),\alpha(t))\,dt\le {\mathcal V}(x)+\eta, \\ \, \\
t_x^\delta(\alpha)\doteq \inf\left\{t>0: \ y_x(t,\alpha)\in\C_\delta\right\}\le  T({\mathcal V}(x)).
\end{array}
\eeq
\end{itemize}}
 
\begin{Theorem}\label{ContV} Assume  {\rm (TPC)}.   If   $\mathcal{V}$ is continuous on $\partial\C$,  then  $\mathcal{V}$ is continuous in its domain and $\lim_{x\to \bar x}\mathcal{V}(x)=+\infty$ for every $\bar x\in\partial  Dom({\mathcal V})$. 
\end{Theorem}

\noindent{\it Proof.}   Owing to Proposition \ref{UCS}, $Dom({\mathcal V})$  is an open subset of $\C^c$ where ${\mathcal V}$ is locally bounded and upper semicontinuous. Let $x_0\in Dom({\mathcal V})$ and let $\nu>0$ be such that $B(x_0,\nu)\subset Dom({\mathcal V})$.   Let $M\doteq \sup \{{\mathcal V}(x): \ \ x\in  B(x_0,\nu)\}+3$ and set $R\doteq {\bf d}(x_0)+\nu$.  Fix $\eta\in]0,1[$. 
By the continuity of ${\mathcal V}$ on the compact set $\partial\C$, there is some $\delta\equiv \delta_{\eta}>0$ such that for any $x\in\C_{2\delta}$ there is a control $\tilde\alpha_x\in{\mathcal A}(x)$ verifying
\bel{etaopt}
\int_0^{t_x(\tilde\alpha_x)}l(y_x(t,\tilde\alpha_x),\tilde\alpha_x(t))\,dt\le \eta.
\eeq
Fix $x_1$, $x_2\in B(x_0,\nu)$    and  assume e.g.  ${\mathcal V}(x_2)\ge{\mathcal V}(x_1)$.  By (TPC), in correspondence of $R$, $\eta$  and $\delta$  introduced above, in view of the definition of $M$, there exist some $T\doteq T(M)$ and  $\alpha_1\in{\mathcal A}(x_1)$  such that 
$$
\int_0^{t_{x_1}(\alpha_1)}l(y_{x_1}(t,\alpha_1),\alpha_1(t))\,dt\le {\mathcal V}(x_1)+\eta<M
$$
and
$$
 {\bf d}(y_{x_1}(\bar t,\alpha_1))<\delta
$$
for some $\bar t\le T\land  t_{x_1}(\alpha_1)$.

\noindent  By standard estimates, all trajectories starting from points   $x\in B(x_0,\nu)$ and corresponding to a control $\alpha$, verify $|y_x(t,\alpha)|\le R'$  \, $\forall t\in[0,T]$ for some   $R'>0$. Moreover, 
$$
|y_{x_2}(t,\alpha)-y_{x_1}(t,\alpha)|\le L|x_2-x_1| \quad \forall t\in[0,T], \quad \forall x_1,\, x_2\in B(x_0,\nu),
$$
for a suitable $L>0$. Therefore,  choosing    $0<\nu'<\nu$ small enough in order to have $ {\bf d}(y_{x_2}(\bar t,\alpha_1))<2\delta$ and denoting by $ \omega(\cdot, R')$ the modulus of $l$ in $B(0,R')$,  for $x_1$, $x_2\in B(x_0,\nu')$, by considering the control $\alpha_2\in {\mathcal A}(x_2)$ given by
$$
\alpha_2(t)=\alpha_1(t)\chi_{[0,\bar t[}+\tilde\alpha_{y_{x_2}(\bar t,\alpha_1)}(t-\bar t)\chi_{[ \bar t, +\infty[}, 
$$
 we get
 $$
\begin{array}{l}
0\le {\mathcal V}(x_2)-{\mathcal V}(x_1)\le \int_0^{\bar t}|l(y_{x_2}(t,\alpha_1),\alpha_1(t))-l(y_{x_1}(t,\alpha_1),\alpha_1(t))|\,dt +3\eta  \\
\le T\, \omega(L|x_2-x_1|, R')+3\eta<4\eta,
\end{array}
$$
 which implies the continuity of  ${\mathcal V}$ by the arbitrariness of $\eta>0$.
 
 To prove that  $\lim_{x\to \bar x \in\partial  Dom({\mathcal V})}\mathcal{V}(x)=+\infty$, assume by contradiction that there are some 
$\bar x\in\partial  Dom({\mathcal V})$,  $M>0$ and $x_n\in Dom({\mathcal V})$ such that  
$|x_n-\bar x|\le 1/n$ and ${\mathcal V}(x_n)\le M$ for any $n\ge 1$. 
  Hence, arguing as in the previous step it not difficult to show that, thanks to the continuity   of ${\mathcal V}$ on the compact set $\partial\C$ and using the (TPC),  starting from some nearly optimal control $\alpha_n\in {\mathcal A}(x_n)$  for $n$ large enough one can construct an admissible control $\alpha \in {\mathcal A}(x)$. Therefore  $\bar x$ does not belong to $\partial  Dom({\mathcal V})$, being $Dom({\mathcal V})$  an open set.


\vv
Hypotheses (SC1), (SC2) introduced in Section \ref{comparison} guarantee the propagation of the continuity of ${\mathcal V}$ from $\partial\C$.

 \begin{Proposition}\label{ContVCS} Assume $\mathcal{V}$  continuous on $\partial\C$ and either  {\rm (SC1)} or {\rm (SC2)}.  Then  $\mathcal{V}$ is continuous in its domain  and $\lim_{x\to \bar x}\mathcal{V}(x)=+\infty$ for every $\bar x\in\partial  Dom({\mathcal V})$.   
\end{Proposition}   

\noindent{\it Proof.}  Fix $R$, $\eta$ and $\delta>0$. Assume first  {\rm (SC1)}.  The UGAS property w.r.t. an invariant set $\C$ yields the existence of a function $\beta\in{\cal KL}$ such that for all   $x\in\C^c$ such that  ${\bf d}(x\le R$ and for  all $\alpha\in{\mathcal A}$,  one has 
$$
{\bf d}(y_x(t,\alpha))\le\beta({\bf d}(x),t)\le\beta(R,t) \qquad \forall t>0
$$
 (see \cite{BaRo}).  Hence   for  all (not  necessarily nearly optimal) controls   $\alpha\in{\mathcal A}$,   the exit  time  $t_x^\delta(\alpha)$ defined as in (TPC) is not greater than $T\doteq\inf\{t>0: \ \  \beta(R,t) \le\delta\}<+\infty$. Therefore, also  in view of Proposition \ref{UCS},  (SC1) together with the continuity of  $\mathcal{V}$   on $\partial\C$  implies   (TPC). 

  If condition (SC2) is assumed, hypothesis  (TPC) is easily fulfilled. Indeed,   for all $x\in\C^c$ and $\alpha\in{\mathcal A}$ with finite cost   $M>0$,    one has  $M=\int_0^{t_x^\delta(\alpha)}l(y_x(t,\alpha),\alpha(t))\,dt\ge c(\delta)\, t_x^\delta(\alpha)$, so that (TPC) is satisfied by choosing   $T(r)\doteq (r+\eta)/c(\delta)$ for any $r\ge0$.

 At this point,  the statement  follows in both cases   from Theorem \ref{ContV}.

   \subsection{Continuity on the target}\label{conttarg}
 Analogously to Proposition \ref{CVT},  by Theorem 1.1 in \cite{MR} it follows that the (LACL)  for $l+\rho$ yields the continuity on    $\partial\C$ of  ${\mathcal V}^\rho_\varepsilon$ for any $\varepsilon\in]0,1]$  and of ${\mathcal U}^\rho$. 
 
In this subsection we provide weaker conditions  sufficient  for the continuity of  ${\mathcal U}^\rho$ on $\partial\C$, involving just a MRF $U$ for $l$.  To this aim  let us  first recall from \cite{MR}  that inequality  (\ref{C1}) is equivalent to (\ref{ECC}) below,  involving a Lagrangian $g$ which, differently from $l$,  is always strictly positive outside the target.  

\begin{Proposition}\label{PMR1} {\rm [MR]}   Let $U$ be a  local  MRF for $h$.  Then $\forall\sigma\in]0,U_0[$   there exists a continuous, strictly increasing  function  $m:]0,\sigma]\to]0,+\infty[$   
 such that, setting
\bel{g}
g(x,a)\doteq {{k}}\, h(x,a)+m(U(x)) \quad  \forall (x,a) \in  \, U^{-1}(]0,\sigma])\times A,
\eeq
 ($k$ the same as in (\ref{C1})),  one has
\bel{ECC}
\min_{ a\in A} \Big\{   \langle p, f(x,a)\rangle +g(x,a) \Big\}\le0  \qquad \forall p\in D^* U(x). 
\eeq  
\end{Proposition}
Following the notations introduced above,  for any $x\in U^{-1}(]0,\sigma])$ let us define for $t\in[0,t_x(\alpha)[$,
$$
c(t)\doteq\int_0^{t}g(y_x(\tau,\alpha),\alpha(\tau))\, d\tau \ \text{ and } \  \bar c\doteq \int_0^{t_x(\alpha)}g(y_x(\tau,\alpha),\alpha(\tau))\, d\tau.
$$
This function is invertible and its inverse, $t(c)$, is a continuous time-change such that
$$
t(c)=\int_0^{c}\frac{dc'}{g(y_x(t(c'),\alpha),\alpha(t(c')))}, \quad c\in[0,\bar c[.
$$
Therefore, $t_x(\alpha)<+\infty$ if and only if
\bel{tfin}
 \int_0^{\bar c}\frac{dc'}{g(y_x(t(c'),\alpha),\alpha(t(c')))}<+\infty
\eeq
and $\int_0^{t_x(\alpha)}\rho(y_x(t,\alpha),\alpha(t))\, dt<+\infty $ if and only if
\bel{rhofin}
 \int_0^{\bar c}\frac{\rho(y_x(t(c'),\alpha),\alpha(t(c')))}{g(y_x(t(c'),\alpha),\alpha(t(c')))}\,dc'<+\infty.
\eeq
  
\noindent In order to state, in Theorem \ref{FTFF} below,  explicit sufficient conditions implying either (\ref{tfin}) or (\ref{rhofin}),  we need the following  definitions.   Owing to  Proposition \ref{PMR1},  fixed a selection  $p(x)\in D^*U(x)$  there is some $a(x)\in A$  such that
 \bel{ce}
 \langle   p(x),\, f(x,a(x))\rangle+ g(x,a(x))\le 0, 
 \eeq
where $g(x,a)= k\,l(x,a)+  m(U(x))$. Let $A(x)$ denote the set of all $a(x)$ verifying (\ref{ce}). 
 For any $s\in]0,\sigma]$,  let us define
\bel{hatg}
\hat g(s)\doteq \inf \{g(x,a): \quad x\in U^{-1}([s,\sigma]), \  a\in A(U^{-1}([s,\sigma]))\},
\eeq
and
\bel{hatrho}
\hat \rho(s)\doteq \sup\{\rho(x,a): \quad x\in U^{-1}(]0,s]), \  a\in A(U^{-1}(]0,s]))\}.
\eeq
Notice that $\hat g(s)\ge   m(s)>0$  and $\hat \rho(s)\le \bar M\doteq\max\{\rho(x,a): \   x\in U^{-1}(]0,\sigma]), \ a\in A\}<+\infty$ for all $s\in]0,\sigma]$. 
As it is not restrictive, let us assume $\hat g$ and $\hat \rho$ continuous. In fact, it is easy to construct  continuous approximations of $\hat g$ and $\hat \rho$ from below and from above, respectively,   verifying all the properties  described above.   
 
\begin{Theorem}\label{FTFF}  Let $\rho$ be a continuous function.  Assume that there exists a   local  MRF $U$ for $l$,  defined in some open set $\Omega$. 
\begin{itemize}
\item[(i)]   If   $\hat g$ and $\hat \rho,$  defined as in (\ref{hatg}), (\ref{hatrho}), respectively, verify  
\bel{CFF}
\int_0^\sigma\frac{ds}{\hat g(s)}=+\infty, \quad    \int_0^\sigma\frac{\hat \rho(s)}{\hat g(s)}\, ds<+\infty,
\eeq
 then ${\mathcal U}^\rho(x)\le  k^{-1}\, U(x)$ \, $\forall  x\in  \Omega$, so that   ${\mathcal U}^\rho$ is continuous on $\partial{\mathcal T}$. Moreover,  ${\mathcal U}^\rho\equiv {\mathcal V}$.

\item[(ii)]
 If  $\hat g$   defined as in (\ref{hatg}) verifies
\bel{CFT}
\int_0^\sigma\frac{ds}{\hat g(s)}<+\infty,
\eeq
then ${\mathcal V}^f(x)\le  k^{-1}\, U(x) \quad \forall  x\in  \Omega$, so that  ${\mathcal V}^f$  is continuous on $\partial{\mathcal T}$.  Moreover,   
${\mathcal V}^f\equiv{\mathcal U}^\rho\equiv{\mathcal V}$ for every $\rho$. 
\end{itemize}
\end{Theorem}

\noindent{\it Proof.}  The results in (i) and (ii)   are  straightforward consequences of Theorem \ref{P11} and of the following Lemma. 

In order to state the lemma, we recall by  Section 3 in \cite{MR} that the existence of a local MRF $U$ for $l$ implies the existence of a  ${\mathcal KL}$ function  $\beta$ such that for any $\eta>0$ there is some $\alpha\in{\mathcal A}$ verifying  
\bel{bbound}
{\bf d}(y_x(t,\alpha))\le \beta({\bf d}(x),t) \qquad \forall t\in\R_+,
\eeq
 and  the estimate 
 \bel{Vprop}
\int_0^{t_x(\alpha)}g(y_x(t,\alpha),\alpha(t))\, dt\le (1+\eta)\,U(x). 
\eeq
 
\begin{Lemma}\label{tecFTFF}  Under the same assumptions and notations of Theorem \ref{FTFF},    if for some $\rho$  condition (\ref{CFF}) holds,   then   $\forall\eta>0$,     $\forall  x\in  U^{-1}(]0,\sigma[)$, there exists  $\alpha\in {\mathcal A}(x)$ verifying  (\ref{bbound}), (\ref{Vprop})
 and also
\bel{th2}
\int_0^{t_x(\alpha)} \rho(y_x(t,\alpha),\alpha(t))\,dt \le (1+\eta)\hat R(U(x)),
\eeq 
where $\hat R(s)\doteq \int_0^s\frac{\hat \rho(s')}{\hat g(s')}\, ds'$  for $s\in]0,\sigma]$.

If  the stronger condition (\ref{CFT}) is satisfied,  then,    $\forall\eta>0$, $\forall  x\in  U^{-1}(]0,\sigma[)$ there exists   $\alpha\in {\mathcal A}(x)$ verifying (\ref{bbound}), (\ref{Vprop}),  and  
\bel{th1}
 t_x(\alpha)\le (1+\eta)\,\hat G(U(x)),
\eeq 
where $\hat G(s)\doteq\int_0^s\frac{ds'}{\hat g(s')}$ for $s\in]0,\sigma]$.
 \end{Lemma}
 
\noindent {\it{Proof.}}
Fix $\eta>0$. Let $(\nu_k)_k\subset]0,1]$ be any infinitesimal, strictly decreasing sequence  such that $\nu_0=1$. Using the same notations of Proposition \ref{FTFF},  for any $x\in U^{-1}(]0,\sigma])$, set $\mu_k\doteq \nu_k\, U(x)$.  In view of the proof of Theorem 1.1 in Section 3 of \cite{MR},  there exist a trajectory-control pair $(y,\alpha):[0,\bar t[\to\C^c\times A'$ and a sequence $t_0\doteq 0<t_1<\dots t_k<t_{k+1}<\dots$,  $\bar t=\lim_{k\to+\infty}t_k$ (possibly, $\bar t=+\infty$),   such that,   for each $k\ge0$, 
$$
U(y(t_k))=\mu_k, \quad U(y(t_{k+1}))<U(y(t))\le U(y(t_k)) \quad\forall t\in[t_k,t_{k+1}[,
$$
$$
\alpha(t)\in A'(U^{-1}(]0,\mu_{k+1}])  \quad\forall t\in[t_k,t_{k+1}[,
$$
 and
\bel{FB}
\int_{t_{k}}^{t_{k+1}}g(y(t),\alpha(t))\,dt\le(1+\eta)[U(y(t_{k}))-U(y(t_{k+1}))].
\eeq
By the definition  (\ref{hatg}) of $\hat g$, 
$$
\int_{t_k}^{t_{k+1}}\hat g(U(y(t_{k+1}))\,dt\le\int_{t_{k}}^{t_{k+1}}g(y(t),\alpha(t))\,dt\le(1+\eta)[\mu_k-\mu_{k+1}]
$$
since $U(y(t_{k}))-U(y(t_{k+1}))=\mu_k-\mu_{k+1}$. Therefore, 
$$
\bar t=\sum_{k=0}^{+\infty}[t_{k+1}-t_k]\le (1+\eta)\sum_{k=0}^{+\infty}\frac{\mu_k-\mu_{k+1}}{\hat g(\mu_{k+1})}.
$$
Assume (\ref{CFT}). Then,  $\forall\varepsilon>0$, since $(\nu_k)_k$ is arbitrary,   we can  choose  $(\mu_k)_k$  so that  
$$\left|\sum_{k=0}^{+\infty}\frac{\mu_k-\mu_{k+1}}{\hat g(\mu_{k+1})}-\int_0^{U(x)}\frac{ds}{\hat g(s)}\right|=\left|\sum_{k=0}^{+\infty}\frac{\mu_k-\mu_{k+1}}{\hat g(\mu_{k+1})}-\hat G(U(x))\right|\le\varepsilon$$  and by the arbitrariness of $\varepsilon,$ $t_x(\alpha)\le\bar t\le (1+\eta) \hat G(U(x))$. The second part of the lemma is thus proved, since  the upper bound  on the cost in (\ref{Vprop}) follows straightforwardly  from (\ref{FB}).

Let now hypothesis (\ref{CFF}) hold. Then 
$$
\int_0^{\bar t}|\alpha(t)|^r\,dt= \sum_{k=0}^{+\infty} \int_{t_{k}}^{t_{k+1}} |\alpha(t)|^r\,dt\le \sum_{k=0}^{+\infty}\hat a^r(U(y(t_{k}))) [t_{k+1}-t_k]
$$
and  by the previous estimates we get
$$
\int_0^{\bar t}|\alpha(t)|^r\,dt\le  \sum_{k=0}^{+\infty}\frac{\hat a^r(\mu_k)}{\hat g(\mu_{k+1})} [\mu_k-\mu_{k+1}].
$$
Since, as it is not restrictive by the arbitrariness of the sequence $(\nu_k)_k$, we can assume that $\hat a(\mu_k)\le C\,\hat a(\mu_{k+1})$  for some $C>0$ for all $k$, the proof can  be easily concluded arguing as above.


\begin{Example}\label{LqEx}{\rm Let  us consider a minimization problem where $U(x)={\bf d}^\gamma(x)$ is a  local MRF in $\C_\sigma$ for some   $\sigma$, $\gamma>0$.  Let us also  assume that, for any $s\in]0,\sigma[$,
$$
\begin{array}{l}
\tilde g(s)\doteq  \inf \{g(x,a): \quad s\le {\bf d}(x)\le\sigma, \  a\in A'({\bf d}^{-1}([s,\sigma]))\}=\bar C_1\,s^{\beta_1}, \\
 \tilde a(s)\doteq \sup\{|a|: \quad a\in A'({\bf d}^{-1}(]0,s]))\}=\bar C_2\,s^{\beta_2},
 \end{array} 
$$
for  some $\bar C_1$, $\bar C_2>0$, and $\beta_1$, $\beta_2\ge0$.  Let us notice that, when $U={\bf d}^\gamma$, using the notations of  Proposition \ref{FTFF}, one has
$$
\hat g(s)=\tilde g(s^{1/\gamma}); \quad \hat a(s)=\tilde a(s^{1/\gamma}).
$$
Then condition (\ref{CFT})  becomes
$$ 
\int_0^\sigma\frac{\,ds}{ s^{ \beta_1/\gamma}}  <+\infty,
$$ 
and it  turns out to be satisfied iff $\beta_1< \gamma$ (without restrictions on the control size, since we can set  $\beta_2=0$).  In case $\beta_1\ge \gamma$, instead, for any $r\ge 1$  the weaker assumption (\ref{CFF}), becoming
$$ 
\int_0^\sigma \frac{\,ds}{s^{( \beta_1-r\,\beta_2)/\gamma}} <+\infty,
$$ 
is verified iff \, $r\,\beta_2>\beta_1-\gamma$. 

Let for instance $\gamma=1$. Then  the stronger condition   $\beta_1<1$  cannot be satisfied  if,  for instance,   $g$ has polynomial growth  around the target  (w.r.t. the distance function) of degree $\beta_1\ge1$.   Notice that $g$ could have $\beta_1<1$ even in cases in which the original Lagrangian $l$ grows as ${\bf d}^\beta(x)$ with $\beta\ge1$.}
 \end{Example}

 \section{Infinite horizon problem}\label{infinite}
Let us introduce the {\it infinite horizon value function,}   defined   for any $x\in\R^n$  as
\begin{equation}\label{calV}
 {\mathcal{ V}}^{\infty}(x)=\inf_{\alpha\in{\mathcal{ A}}}{\mathcal J}(+\infty,x,\alpha).
\end{equation}
Clearly, for any target $\C$ one has  $\mathcal{ V}^m\le {\mathcal{ V}}^{\infty}$  but, choosing   $\C$ in a suitable way,  ${\mathcal{ V}}^{\infty}$ does actually coincide with  $\mathcal{ V}^m$.  As a consequence, in this case   all the results obtained   for ${\mathcal V}$ yield analogous results for ${\mathcal V}^\infty$ as soon as we have ${\mathcal V}\equiv {\mathcal V}^m$. We recall by Proposition \ref{CorUn}  that either hypothesis (SC1) or (SC2) implies such an equality. 
 For instance, we have
\begin{Proposition}\label{Vinfty} 
   If  $\C\times\{0\}$ is a viability set \footnote{Let $F(x)\doteq \{(f(x,a),l(x,a)): \ a\in    A\}$.     Any closed subset $K\subset\R^n\times\R$ will be called   a {\it viability set} for $(f,l)$ if  for any $(x_0,\lambda_0)\in K$ there is a solution $(y,\lambda)$ of 
the differential inclusion
$$
(\dot y(t), \dot\lambda(t))\in F(y(t))\qquad   t\ge 0
$$
 such that  $(y(0),\lambda(0))=(x_0,\lambda_0)$ and $(y(t),\lambda(t))\in K$ \ $\forall t>0$ (see e.g. \cite{AF}).}  for $(f,l)$,
then $\mathcal{V}^{\infty} \equiv{\mathcal{ V}}^m$.
 \end{Proposition}
As an immediate consequence we get the following
\begin{Corollary}\label{Vinftyhp}  If  $\C\times\{0\}$ is a viability set for $(f,l)$ and either {\rm (SC1)} or {\rm (SC2)} holds true for $\C$, then  $\mathcal{V}^{\infty} \equiv{\mathcal{ V}}^m\equiv \mathcal{V}$.
\end{Corollary}

 \noindent{\it Proof of Proposition \ref{Vinfty}.}   
 Let $x\in\R^n$ be such that ${\mathcal{ V}}^m(x)<+\infty$ (otherwise, the equality $\mathcal{V}^{\infty}(x)={\mathcal{ V}}^m(x)=+\infty$ is trivial). For any $\eta>0$ there is some $\hat\alpha\in{\mathcal A}$ such that
\bel{infeq}
 \int_0^{t_x(\hat\alpha)} l(y_{x}(t,\hat\alpha),\hat\alpha(t))\,dt\le {\mathcal{ V}}^m(x) +\eta.
\eeq
If $t_x(\hat\alpha)=+\infty$,   (\ref{infeq}) implies that ${\mathcal{ V}}^{\infty}(x)\le {\mathcal{ V}}^m(x)+\eta$.  If $t_x(\hat\alpha)<+\infty$, set $\bar x\doteq y_{x}(t_x(\hat\alpha),\hat\alpha)$.
 Since  $\C\times\{0\}$ is viable for $(f,l)$ and $\bar x\in\C$,  there exists a control $\bar\alpha\in{\mathcal A}$ such  that 
$$
y_{\bar x}(t,\bar\alpha)\in\C \quad\text{and}\quad \int_0^t l(y_{\bar x}(t,\bar\alpha),\bar\alpha(t))\,dt=0 \quad \forall t>0
 $$
(incidentally, this argument yields that $\mathcal{V}^{\infty}\equiv0$ in $\C$).
Therefore considering  the control $\alpha(t)\doteq \hat\alpha\chi_{[0,t_x(\hat\alpha)[}(t)+\bar\alpha(t-t_x(\hat\alpha))\chi_{[t_x(\hat\alpha),+\infty[}(t)$ for all $t\ge0$, we get that  $\int_0^{+\infty} l(y_{x}(t,\alpha),\alpha(t))\,dt=\int_0^{t_x(\hat\alpha)} l(y_{x}(t,\hat\alpha),\hat\alpha(t))\,dt$. Hence   (\ref{infeq}) implies  again that ${\mathcal{ V}}^{\infty}(x)\le {\mathcal{ V}}^m(x)+\eta$. By the arbitrariness of $\eta>0$,  this yields the equality  $\mathcal{V}^{\infty} \equiv{\mathcal{ V}}^m$ and the proof is concluded.    

\begin{Remark}\label{Rviab}{\rm   Viability sufficient  conditions can be found e.g.  in \cite{AF}.  In particular, we recall that, 
if   ${\mathcal T}$ and the sets   $F(x)\doteq \{(f(x,a),l(x,a)): \ a\in    A\}$ $\forall x\in{\mathcal T}_\sigma$ are convex and the (LACL) for $l$ holds,  the viability  hypothesis in Proposition \ref{Vinfty} is satisfied by well known results of convex analysis.  Actually in this case there exists an equilibrium point $\bar x\in{\mathcal T}$ for   system (\ref{S}).   For nonconvex  sets,  the implication is in general false. }
\end{Remark}
  
\begin{Remark}{\rm  The viability hypothesis in Proposition \ref{Vinfty}  is a sufficient condition in order to have  $\mathcal{V}^{\infty}\equiv0$ on $\C$.   Actually, by the proof above one easily deduces that,  if ${\mathcal Z}_\infty\doteq\{x\in\R^n: \ \   {\mathcal{ V}}^{\infty}(x)=0\}\ne\emptyset$ then
$$
\text{ {\it $\mathcal{V}^{\infty} \equiv{\mathcal{ V}}^m$ for any  target   $\C\subset{\mathcal Z}_\infty$.}}
$$
}
\end{Remark}

  Under the hypotheses of Corollary \ref{Vinftyhp}, all the results about uniqueness, stability and continuity obtained in the previous sections for ${\mathcal V}$ give rise to analogous results for the infinite horizon value function, $\mathcal{V}^{\infty}$, which we omit for the sake of brevity.   For instance,  in view of Remark \ref{Rviab},   Corollary \ref{Unexpl} implies the following.
   \begin{Corollary}\label{VinftyUn}  Let either $\C\times\{0\}$ be a viability set for $(f,l)$ or  ${\mathcal T}$ and the sets   $F(x)\doteq \{(f(x,a),l(x,a)): \ a\in    A\}$ $\forall x\in{\mathcal T}_\sigma$ be convex.  Assume  {\rm (CV)},  {\rm (LACL)}  for $l$ and either {\rm (SC1)} or {\rm (SC2)}. Then 
   \begin{itemize}
   \item[(i)]  there is a unique nonnegative viscosity solution $W$  to  {\rm(BVP$\mathcal{K}$)}, which turns out to be continuous. Moreover, ${\mathcal V}^\infty\equiv\Psi^{-1}(W)=-\log(1-W)$ and $Dom(\mathcal{V}^\infty)=\{x:\ \ W(x)<1\}$;
    \item[(ii)]    $(\mathcal{V}^\infty, Dom(\mathcal{V}^\infty))$  is  the unique nonnegative  viscosity solution to  {\rm({BVP})}  among the pairs $(u,\Omega)$, where $u$ satisfies {\rm (\ref{BCn})}. 
Moreover,  $\mathcal{V}^\infty$ is continuous.
    \end{itemize}  
 If  {\rm (CV)} is not satisfied,  the above uniqueness results hold just among the continuous functions.     
 \end{Corollary}
 
We conclude by observing that Examples \ref{Un1} and \ref{Un3} are in fact examples of infinite horizon problems.
The value function ${\mathcal V}$ coincides with ${\mathcal V}^{\infty}$ and it is the unique nonnegative solution of the (BVP) with   $\C=\{0\}$. 

\begin{Remark} \label{infH} {\rm 
We can cover affine-quadratic  problems, where
 $$
 l(x,a)=x^T Q\, x+a^T R \,a, \qquad f(x,a)=A(x)+\langle B(x),a\rangle, \quad x\in\R^n, a\in A\subset \R^m
 $$
 with $Q$ and $R$ symmetric and {\it positive semi-definite} matrices and the control set $A$ convex. In the literature, $Q$ is usually assumed to be positive definite. In this case, (SC2) for $\C=\{0\}$ is trivially satisfied,  so that  one has ${\mathcal V}\equiv {\mathcal V}^m$.

 For such problems we have uniqueness among the nonnegative solutions of (BVP) for $\C=\{0\}$  assuming the (LACL)  and either $Q$ is positive definite or $f$ has the UGAS property w.r.t. $\{0\}$, namely (SC1) holds.
 
Such kind of problems are widely studied, both for unbounded and bounded controls. In the last  case  they are also known in the literature  as constrained or saturated  problems (see e.g. \cite{G}).  

}
\end{Remark}

 \end{large}
\end{document}